\documentclass[12pt,a4paper,fleqn]{article}
\usepackage{a4wide,amsfonts,amsmath,latexsym,amssymb,euscript,graphicx,units,mathrsfs}

\usepackage{graphicx}
\usepackage{color}
\usepackage{amssymb}
\usepackage{amssymb}
\usepackage[T1]{fontenc}
\usepackage{latexsym}
\usepackage{xypic}
\usepackage{eufrak}
\usepackage{euscript}
\usepackage{amsfonts,amsmath}
\usepackage{verbatim}
\usepackage{fancyhdr}
\usepackage[english]{babel}
\usepackage{mathrsfs}
\usepackage{units}

\newtheorem{prop}{Proposition}[section]

\newtheorem{lemme}[prop]{Lemma}
\newtheorem{rem}[prop]{Remark}
\newtheorem{thm}[prop]{Theorem}
\newtheorem{defi}[prop]{Definition}

\renewcommand{\geq}{\geqslant}
\def\leq{\leqslant}

\newcommand{\R}{\mathbb{R}}

\def\1{{\mathbf{1}}}

\def\1{{\mathbf{1}}}
\def\0.5{{\frac{1}{2}}}

\newcommand{\fin}
{ \vspace{-0.6cm}
\begin{flushright}
\mbox{$\Box$}
\end{flushright}
\noindent }

\newcommand{\qed}{\nopagebreak\hspace*{\fill}
{\vrule width6pt height6ptdepth0pt}\par}

\begin{document}

\begin{center}
{\large{\bf
Yet another proof of the Nualart-Peccati criterion
}}\\~\\
by Ivan Nourdin\footnote{Institut \'Elie Cartan, Universit\'e Henri Poincar\'e, BP 70239, 54506 Vandoeuvre-l\`es-Nancy, France,
{\tt inourdin@gmail.com}}\footnote{supported in part by the two following (french) ANR grants: `Exploration des Chemins Rugueux'
[ANR-09-BLAN-0114] and `Malliavin, Stein and Stochastic Equations with Irregular Coefficients'
[ANR-10-BLAN-0121].}\\
{\it Universit\'e Nancy 1}\\~\\
\small{{\it This version}: December 15th, 2011}\\
\end{center}

{\small
\noindent
{\bf Abstract:}
In \cite{nualart-peccati}, Nualart and Peccati showed that, surprisingly, the convergence in distribution
of a normalized sequence of multiple Wiener-It\^o integrals towards a standard Gaussian law
is equivalent to
convergence of just the fourth moment to 3.
In \cite{withspeicher}, this result is extended to a sequence of multiple Wigner integrals, in the context
of free Brownian motion.
The goal of the present paper is to offer an elementary, unifying proof of these two results.
The only advanced, needed tool
is the product formula for multiple integrals. Apart from this formula,
the rest of the proof only relies on soft combinatorial arguments.\\


\section{Introduction}

The following surprising result, proved in \cite{nualart-peccati}, shows that the convergence in distribution
of a normalized sequence of multiple Wiener-It\^o integrals towards a standard Gaussian law
is equivalent to
convergence of just the fourth moment to 3.
\begin{thm}
[Nualart-Peccati]\label{nunugio}
Fix an integer $p\geq 2$. Let $\{B(t)\}_{t\in [0,T]}$ be a classical Brownian motion,
and let $(F_n)_{n\geq 1}$
be a sequence of multiple integrals of
the form
\begin{equation}\label{fn}
F_n=\int_{[0,T]^p} f_n(t_1,\ldots,t_p)dB(t_1)\ldots dB(t_p),
\end{equation}
where each $f_n\in L^2([0,T]^p;\R)$ is symmetric (it is not a restrictive assumption). Suppose moreover that $E[F_n^2]\to 1$ as $n\to\infty$.
Then, as $n\to\infty$, the following two assertions are equivalent:
\begin{enumerate}
\item[(i)] The sequence $(F_n)$ converges in distribution to $B(1)\sim N(0,1)$;
\item[(ii)] $E[F_n^4]\to E[B(1)^4]=3$.
\end{enumerate}
\end{thm}
In \cite{nualart-peccati}, the original proof of $(ii)\Rightarrow (i)$ relies on tools from
Brownian stochastic analysis. Precisely, using the symmetry of $f_n$, one can rewrite $F_n$ as
\[
F_n=p!\int_0^T dB(t_1)\int_0^{t_1} dB(t_2)\ldots \int_0^{t_{p-1}}dB(t_p) f_n(t_1,\ldots,t_p),
\]
and then make use of the Dambis-Dubins-Schwarz theorem to transform it into
$F_n = \beta^{(n)}_{\langle F_n\rangle}$, where $\beta^{(n)}$ is a classical Brownian motion and
\begin{equation}\label{dds}
\langle F_n\rangle = p!^2\int_0^T dt_1 \left(\int_0^{t_1} dB(t_2)\ldots \int_0^{t_{p-1}}dB(t_p)
f_n(t_1,\ldots,t_p)\right)^2.
\end{equation}
Therefore, to get that $(i)$ holds true, it is now enough to prove that $(ii)$ implies $\langle F_n\rangle \overset{L^2}{\to} 1$,
which is exactly what Nualart and Peccati did in \cite{nualart-peccati}.
\\

Since the publication of \cite{nualart-peccati}, several researchers have been
interested in understanding more deeply why
Theorem \ref{nunugio} holds. Let us mention some works in this direction:\\

1. In \cite{nunuol}, Nualart and Ortiz-Latorre gave another proof of Theorem \ref{nunugio}
using exclusively the tools of Malliavin calculus. The main ingredient of their proof is
the identity $\delta D= -L$, where $\delta$, $D$ and $L$ are basic operators in Malliavin calculus.\\

2. Based on the ideas developed in \cite{stein-PTRF},
the following bound is shown in \cite[Theorem 3.6]{stein-survey} (see also
\cite{NPR}): if $E[F_n^2]=1$,
then
\begin{equation}\label{nou-pec}
\sup_{A\in\mathscr{B}(\R)}\left|P[F_n\in A]-\frac{1}{\sqrt{2\pi}}\int_{A} e^{-u^2/2}du\right|\leq 2\sqrt{\frac{p-1}{3p}}\sqrt{|E[F_n^4]-3|}.
\end{equation}
Of course, with (\ref{nou-pec}) in hand, it is totally straightforward to obtain
Theorem \ref{nunugio} as a corollary.
However, the proof of (\ref{nou-pec}), albeit not that  difficult, requires the knowledge of both Malliavin calculus
and Stein's method.\\

3. By using the tools of Malliavin calculus, Peccati and I
computed in \cite{cumulants} a new expression for the cumulants of $F_n$, in terms of the contractions
of the kernels $f_n$. As an immediate byproduct of this formula, we are able to
recover Theorem \ref{nunugio}, see \cite[Theorem 5.8]{cumulants}
for the details.
See also \cite{nono} for an extension in the multivariate setting.\\

4. In \cite{noncentral}, Theorem \ref{nunugio} is extended to
the case where, instead of $B(1)\sim N(0,1)$ in the limit,
a centered chi-square random variable, say $Z$, is considered.
More precisely, it is proved in this latter reference that an adequably normalized sequence $F_n$ of the form (\ref{fn})
converges  in distribution towards $Z$
 if and only if $E[F_n^4]-12E[F_n^3]\to E[Z^4]-12E[Z^3]$. Here again, the proof is based on
the use of the basic operators of Malliavin calculus.\\

5. The following result, proved in \cite{withspeicher}, is the exact analogue of Theorem \ref{nunugio},
but
in
the situation where the classical Brownian motion
$B$ is replaced by its {\it free} counterpart $S$.
\begin{thm}
[Kemp-Nourdin-Peccati-Speicher]\label{knps}
Fix an integer $p\geq 2$. Let $\{S(t)\}_{t\in [0,T]}$ be a free Brownian motion,
and let $(F_n)_{n\geq 1}$
be a sequence of multiple integrals of
the form
\[
F_n=\int_{[0,T]^p} f_n(t_1,\ldots,t_p)dS(t_1)\ldots dS(t_p),
\]
where each $f_n\in L^2([0,T]^p;\R)$ is mirror symmetric (that is, satisfies
$f_n(t_1,\ldots,t_p)=f_n(t_p,\ldots,t_1)$ for all $t_1,\ldots,t_p\in[0,1]$).
Suppose moreover that $E[F_n^2]\to 1$ as $n\to\infty$.
Then, as $n\to\infty$, the following two assertions are equivalent:
\begin{enumerate}
\item[(i)] For all $k\geq 3$, $E[F_n^k]\to E[S(1)^k]$;
\item[(ii)] $E[F_n^4]\to E[S(1)^4]=2$.
\end{enumerate}
\end{thm}
The proof of Theorem \ref{knps} contained in \cite{withspeicher}
is based on the use of combinatorial features related to the free probability realm,
including non-crossing pairing and partitions. \\

Thus, there is already several proofs of Theorem \ref{nunugio}. Each of them
has its own interest, because it
allows to understand more deeply a particular aspect of this beautiful result.
On the other hand, all these proofs require at some point to deal with
sophisticated tools, such as stochastic Brownian analysis, Malliavin calculus or Stein's method.\\

The goal of this paper is to offer an elementary, unifying proof of both
Theorems \ref{nunugio} and \ref{knps}.
As anticipated, the only advanced result we will need
is the {\it product formula} for multiple integrals,
that is, the explicit expression for the product of two multiples integrals
of order $p$ and $q$, say, as a linear combination of multiple integrals of
order less or equal to $p+q$. Apart from this formula,
the rest of the proof only relies on `soft' combinatorial arguments.\\

The level of our paper is (hopefully) available to any good student.
From our opinion however, its interest is not only to provide a new, simple proof of a known
result. It is indeed noteworthy that the number
of required tools has been reduced
to its maximum (the product formula being essentially the only one we need),
so that our approach might represent a valuable
strategy to follow in order to generalize Theorem(s) \ref{nunugio} (and \ref{knps})
in other situations.
For instance, let us mention that the two works \cite{poisson,tetilla} have indeed
followed our line of reasoning, and successfully extended Theorem \ref{knps} in the case where the limit is
the free Poisson distribution and the (so-called) tetilla law respectively.\\

The rest of the paper is organized as follows. Section 2 deals with some preliminary results.
Section 3 contains our proof of Theorem \ref{knps}, whereas Section 4 is devoted
to the proof of Theorem \ref{nunugio}.

\section{Preliminaries}

\subsection{Multiple integrals with respect to classical Brownian motion}
In this section, our main reference is Nualart's book \cite{nualart-book}. To simplify the exposition,
without loss of generality we fix the time horizon to be $T=1$.

Let $\{B(t)\}_{t\in[0,1]}$ be a classical Brownian motion,
that is, a stochastic process defined on a probability space $(\Omega,\mathscr{F},P)$,
starting from 0, with independent increments,
and such that $B(t)-B(s)$ is
a centered Gaussian random variable with variance $t-s$ for all $t\geq s$.

For a given real-valued kernel $f$ belonging to $L^2([0,1]^p)$,
let us quickly sketch out the construction of the {\it multiple Wiener-It\^o integral} of $f$ with respect to $B$,
written
\begin{equation}\label{multq}
I_p(f)=\int_{[0,1]^p} f(t_1,\ldots,t_p)dB(t_1)\ldots dB(t_p)
\end{equation}
in the sequel.
(For the full details, we refer the reader to the classical reference \cite{nualart-book}.)
Let $D^p\subset[0,1]^p$ be the collection of all diagonals, i.e.
\begin{equation}\label{dp}
D^p=\{(t_1,\ldots,t_p)\in[0,1]^p:\,t_i=t_j\mbox{ for some $i\neq j$}\}.
\end{equation}
As a first step, when $f$ has the form of a characteristic function $f={\bf 1}_A$, with
$A=[u_1,v_1]\times\ldots\times [u_p,v_p]\subset[0,1]^p$
such that $A\cap D^p=\emptyset$, the $p$th multiple integral of $f$
is defined
by
\[
I_p(f)=(B(v_1)-B(u_1))\ldots (B(v_p)-B(u_p)).
\]
Then, this definition is extended by linearity to simple functions of the form
$
f=\sum_{i=1}^k \alpha_i {\bf 1}_{A_i},
$
where
$
A_i=[u^i_1,v^i_1]\times\ldots\times [u^i_p,v^i_p]
$
are disjoint $p$-dimensional rectangles as above which do not meet the diagonals.
Simple computations show that
\begin{eqnarray}
E[I_p(f)]&=&0\label{isom1}\\
I_p(f)&=&I_p(\widetilde{f})\label{isom2}\\
\label{isom}
E[I_p(g)I_p(f)]&=&p!\langle \widetilde{g},\widetilde{f}\rangle_{L^2([0,1]^p)}.
\end{eqnarray}
Here, $\widetilde{f}\in L^2([0,1]^p)$ denotes the symmetrization of $f$, that is, the symmetric function
canonically associated to $f$, given by
\begin{equation}\label{norm}
\widetilde{f}(t_1,\ldots,t_p)=\frac{1}{p!}\sum_{\pi\in\mathfrak{S}_p} f(t_{\pi(1)},\ldots,t_{\pi(p)}).
\end{equation}
Since each $f\in L^2([0,1]^p)$ can be approximated in $L^2$-norm by simple functions,
we can finally extend the definition of (\ref{multq}) to all $f\in L^2([0,1]^p)$. Note that, by construction,
(\ref{isom1})-(\ref{isom}) is still true in this general setting.
Then, one easily sees that, in addition,
\begin{equation}\label{isom-dif}
E[I_p(f)I_q(g)]=0\,\mbox{ for any $p\neq q$, $f\in L^2([0,1]^p)$ and $g\in L^2([0,1]^q)$}.
\end{equation}

Before being in position to state the product formula for two multiple integrals, we need to introduce the following quantity.
\begin{defi}\label{contr}
For symmetric functions $f\in L^2([0,1]^p)$ and $g\in L^2([0,1]^q)$, the contractions
\[
f\otimes_r g\in L^2([0,1]^{p+q-2r})\quad\mbox{($0\leq r\leq \min(p,q)$)}
\]
are the (not necessarily symmetric) functions given by
\begin{eqnarray*}
&&f\otimes_r g(t_1,\ldots,t_{p+q-2r}):=\\
&& \hskip2cm\int_{[0,1]^{r}} f(t_1,\ldots,t_{p-r},s_1,\ldots,s_r)g(t_{p-r+1},\ldots,t_{p+q-2r},s_1,\ldots,s_r)
ds_1\ldots ds_r.
\end{eqnarray*}
By convention, we set $f\otimes_0 g= f\otimes g$, the tensor product of $f$ and $g$.
\end{defi}
The symmetrization of $f\otimes_r g$ is written $f\widetilde{\otimes}_r\, g$.
Observe that
$f\otimes_p g=f\widetilde{\otimes}_p\, g=\langle f,g\rangle_{L^2([0,1]^p)}$ whenever $p=q$.
Also, using Cauchy-Schwarz inequality, it is immediate to prove that
\[
\|f\otimes_r g\|_{L^2([0,1]^{p+q-2r})}\leq \|f\|_{L^2([0,1]^p)}\|g\|_{L^2([0,1]^q)}
\]
for all $r=0,\ldots,\min(p,q)$. (It is actually an equality for $r=0$.)
Moreover, a simple application of the triangle inequality leads to
\[
\|f \widetilde{\otimes}_r\, g\|_{L^2([0,1]^{p+q-2r})}\leq \|f\otimes_r g\|_{L^2([0,1]^{p+q-2r})}.
\]

We can now state the {\it product formula}, which is the main ingredient of our proof of Theorem
\ref{nunugio}.
By taking the expectation in (\ref{product}), observe that we recover both (\ref{isom}) and (\ref{isom-dif}).
\begin{thm}\label{product-thm}
For symmetric functions $f\in L^2([0,1]^p)$ and $g\in L^2([0,1]^q)$, we have
\begin{equation}\label{product}
I_p(f)I_q(g)=\sum_{r=0}^{\min(p,q)} r!\binom{p}{r}\binom{q}{r} I_{p+q-2r}(f\widetilde{\otimes}_r\,g).
\end{equation}
\end{thm}

\subsection{Multiple integrals with respect to free Brownian motion}
In this section, our main references are: $(i)$
the monograph \cite{nica-speicher} by Nica and Speicher for the generalities about
free probability; $(ii)$
the paper \cite{BiaSpeich-PTRF} by Biane and Speicher for the free stochastic analysis.
We refer the reader to them for any unexplained notion or result.

Let $\{S(t)\}_{t\in[0,1]}$ be a free Brownian motion,
that is, a stochastic process defined on a non-commutative probability space $(\mathscr{A},E)$,
starting from 0,
with freely independent increments, and such that $S(t)-S(s)$ is
a centered semicircular random variable with variance $t-s$ for all $t\geq s$.
We may think of free Brownian motion as `infinite-dimensional matrix-valued Brownian motion'.
For more details about the construction and features of $S$, see \cite[Section 1.1]{BiaSpeich-PTRF}
and the references therein.

When $f\in L^2([0,1]^p)$ is real-valued, we write $f^*$ to indicate the function of $L^2([0,1]^p)$ given by
$f^*(t_1,\ldots,t_p)=f(t_p,\ldots,t_1)$. (Hence, to say that $f_n$ is mirror-symmetric
in Theorem \ref{knps} means that $f_n=f_n^*$.)
We quickly sketch out the construction of the {\it multiple Wigner
integral} of $f$ with respect to  $S$.
Let $D^p\subset[0,1]^p$ be the collection of all diagonals, see (\ref{dp}).
For a characteristic function $f={\bf 1}_A$, where $A\subset[0,1]^p$ has the form
$
A=[u_1,v_1]\times\ldots\times [u_p,v_p]
$
with $A\cap D^p=\emptyset$, the $p$th multiple Wigner integral of $f$, written
\[
I_p(f)=\int_{[0,1]^p} f(t_1,\ldots,t_p)dS(t_1)\ldots dS(t_p),
\]
is defined
by
\[
I_p(f)=(S(v_1)-S(u_1))\ldots (S(v_p)-S(u_p)).
\]
Then, as in the previous section we extend this definition by linearity to simple functions of the form
$
f=\sum_{i=1}^k \alpha_i {\bf 1}_{A_i},
$
where
$
A_i=[u^i_1,v^i_1]\times\ldots\times [u^i_p,v^i_p]
$
are disjoint $p$-dimensional rectangles as above which do not meet the diagonals.
Simple computations show that
\begin{eqnarray}
E[I_p(f)]&=&0\label{isomfree1}\\
E[I_p(f)I_p(g)]&=&\langle f,g^*\rangle_{L^2([0,1]^p)}.\label{isomfree}
\end{eqnarray}
By approximation, the definition of $I_p(f)$ is extended to all $f\in L^2([0,1]^p)$,
and (\ref{isomfree1})-(\ref{isomfree}) continue to hold true in this more general setting.
It turns out that
\begin{equation}\label{isom-dif-free}
E[I_p(f)I_q(g)]=0\,\mbox{ for $p\neq q$, $f\in L^2([0,1]^p)$ and $g\in L^2([0,1]^q)$}.
\end{equation}

Before giving the product formula in the free context, we need to introduce the analogue
for Definition \ref{contr}.
\begin{defi}\label{contr-free}
For functions $f\in L^2([0,1]^p)$ and $g\in L^2([0,1]^q)$, the contractions
\[
f\overset{r}{\frown} g\in L^2([0,1]^{p+q-2r})\quad\mbox{($0\leq r\leq \min(p,q)$)}
\]
are the functions given by
\begin{eqnarray*}
&&f\overset{r}{\frown} g(t_1,\ldots,t_{p+q-2r}):=\\
&& \hskip2cm\int_{[0,1]^{r}} f(t_1,\ldots,t_{p-r},s_1,\ldots,s_r)g(s_r,\ldots,s_1,t_{p-r+1},\ldots,t_{p+q-2r})
ds_1\ldots ds_r.
\end{eqnarray*}
By convention, we set $f\overset{0}{\frown} g= f\otimes g$, the tensor product of $f$ and $g$.
\end{defi}
Observe that
$f\overset{p}{\frown} g=\langle f,g^*\rangle_{L^2([0,1]^p)}$ whenever $p=q$.
Also, using Cauchy-Schwarz, it is immediate to prove that
$\|f\overset{r}{\frown} g\|_{L^2([0,1]^{p+q-2r})}\leq \|f\|_{L^2([0,1]^p)}\|g\|_{L^2([0,1]^q)}$
for all $r=0,\ldots,\min(p,q)$. (It is actually an equality for $r=0$.)

We can now state the {\it product formula}
in the free context, which turns out to be simpler compared to the classical case (Theorem \ref{product-thm}).
\begin{thm}
For functions $f\in L^2([0,1]^p)$ and $g\in L^2([0,1]^q)$, we have
\begin{equation}\label{prodfree}
I_p(f)I_q(g)=\sum_{r=0}^{\min(p,q)} I_{p+q-2r}(f\overset{r}{\frown}g).
\end{equation}
\end{thm}

\section{Proof of Theorem \ref{knps}}
Let the notation and assumptions of Theorem \ref{knps} prevail.
Without loss of generality, we may assume that $E[F_n^2]=1$ for all $n$
(instead of $E[F_n^2]\to 1$ as $n\to\infty$).
Moreover, because $f_n=f_n^*$, observe that
$\|f_n\|^2_{L^2([0,1]^p)}=E[F_n^2]=1$.

It is trivial that $(i)$ implies $(ii)$.
Conversely, assume that $(ii)$ is in order, and let us prove that $(i)$ holds.
Fix an integer $k\geq 3$. Iterative applications of the product formula (\ref{prodfree}) leads to
\begin{equation}\label{kthmomentwithoutexpectation}
F_n^k=I_p(f_n)^k=\sum_{(r_1,\ldots,r_{k-1})\in A_k}
I_{kp-2r_1-\ldots-2r_{k-1}}\big(
f_n\overset{r_1}{\frown}\ldots\overset{r_{k-1}}{\frown}f_n
\big),
\end{equation}
where
\begin{eqnarray*}
A_k&=&\big\{
(r_1,\ldots,r_{k-1})\in \{0,1,\ldots,p\}^{k-1}:\,r_2\leq 2p-2r_1,\,\,r_3\leq 3p-2r_1-2r_2,\ldots,\\
&&\hskip7cm r_{k-1}\leq(k-1)p-2r_1-\ldots-2r_{k-2}
\big\}.
\end{eqnarray*}
In order to simplify the exposition, note that we have removed the brackets in the writing of
$f_n
\overset{r_1}{\frown}
\ldots
\overset{r_{k-1}}{\frown}
f_n$. We use the implicit convention
that
these quantities are always  defined iteratively from the left to the right. For instance,
$f_n\overset{r_1}{\frown}f_n \overset{r_{2}}{\frown}f_n \overset{r_{3}}{\frown}f_n$ actually
stands for $((f_n\overset{r_1}{\frown}f_n) \overset{r_{2}}{\frown}f_n) \overset{r_{3}}{\frown}f_n$.

By taking the expectation in (\ref{kthmomentwithoutexpectation}), we deduce that
\begin{equation}\label{kthmoment}
E[F_n^k]=\sum_{(r_1,\ldots,r_{k-1})\in B_k}
f_n\overset{r_1}{\frown}\ldots\overset{r_{k-1}}{\frown}f_n,
\end{equation}
with $B_k=\big\{(r_1,\ldots,r_{k-1})\in A_k:\,2r_1+\ldots+2r_{k-1}=kp\big\}$.
We decompose $B_k$ as $C_k\cup E_k$, with $C_k=B_k\cap\{0,p\}^{k-1}$ and $E_k=B_k\setminus C_k$.
We then have, for all $k\geq 3$,
\begin{equation}
E[F_n^k]= \sum_{(r_1,\ldots,r_{k-1})\in C_k}
f_n\overset{r_1}{\frown}  \ldots\overset{r_{k-1}}{\frown}  f_n+\sum_{(r_1,\ldots,r_{k-1})\in E_k}
f_n\overset{r_1}{\frown}  \ldots\overset{r_{k-1}}{\frown}  f_n.\label{doublesum}
\end{equation}
Lemmas \ref{lm2} and \ref{lm3} imply together that
the first sum in (\ref{doublesum}) is equal to $E[S(1)^k]$.
Moreover, by Lemma \ref{lm1} and because $(ii)$ is in order, we have
that $\|f_n\overset{r}{\frown} f_n\|_{L^2([0,1]^{2p-2r})}\to 0$ for all $r=1,\ldots,p-1$.
Hence,
the second sum in (\ref{doublesum}) must converge to zero by Lemma \ref{lm4}.
Thus, $(i)$ is in order, and the proof of the theorem is concluded.\fin

\begin{lemme}\label{lm1}
We have $E[F_n^4]=2+\sum_{r=1}^{p-1} \|f_n \overset{r}{\frown} f_n\|^2_{L^2([0,1]^{2p-2r})}$.
\end{lemme}
{\it Proof}. The product formula (\ref{prodfree}) yields
$
F_n^2= \sum_{r=0}^p I_{2p-2r}(f_n\overset{r}{\frown} f_n).
$
Using (\ref{isomfree})-(\ref{isom-dif-free}), we infer
\begin{eqnarray*}
E[F_n^4] &=& \|f_n\otimes f_n\|^2_{L^2([0,1]^{2p})} + \big(\|f_n\|^2_{L^2([0,1]^p)}\big)^2+
\sum_{r=1}^{p-1}
\langle f_n\overset{r}{\frown} f_n, (f_n\overset{r}{\frown} f_n)^*\rangle_{L^2([0,1]^{2p-2r})}
  \notag\\
&=&  2\|f_n\|^4_{L^2([0,1]^p)}+\sum_{r=1}^{p-1}\|f_n\overset{r}{\frown} f_n\|^2_{L^2([0,1]^{2p-2r})}
=2+\sum_{r=1}^{p-1}\|f_n\overset{r}{\frown} f_n\|^2_{L^2([0,1]^{2p-2r})}, \label{fourthpower}
\end{eqnarray*}
since $\|f_n\|^2_{L^2([0,1]^p)}=1$ and
\begin{eqnarray*}
&&f_n\overset{r}{\frown} f_n(t_1,\ldots,t_{2p-2r})\\
&=&\int_{[0,1]^{r}} f_n(t_1,\ldots,t_{p-r},s_1,\ldots,s_r)f_n(s_r,\ldots,s_1,t_{p-r+1},\ldots,t_{2p-2r})ds_1\ldots ds_r
\\
&=&\int_{[0,1]^{r}} f_n(s_r,\ldots,s_1,t_{p-r},\ldots,t_1)f_n(t_{2p-2r},\ldots,t_{p-r+1},s_1,\ldots,s_r)ds_1\ldots ds_r
\\
&=&f_n\overset{r}{\frown} f_n(t_{2p-2r},\ldots,t_1)=(f_n\overset{r}{\frown} f_n)^*(t_1,\ldots,t_{2p-2r}).
\end{eqnarray*}
\fin
\begin{lemme}\label{lm2}
For all $k\geq 3$, the cardinality of $C_k$ coincides with $E[S(1)^k]$.
\end{lemme}
{\it Proof}. By dividing all the $r_i$'s by $p$, one get that
\begin{eqnarray*}
C_k&\overset{\rm bij.}{\equiv}&\widetilde{C}_k:=\big\{
(r_1,\ldots,r_{k-1})\in \{0,1\}^{k-1}:\,r_2\leq 2-2r_1,\,\,r_3\leq 3-2r_1-2r_2,\ldots,\\
&&\hskip3.5cm r_{k-1}\leq k-1-2r_1-\ldots-2r_{k-2},\,2r_1+\ldots+2r_{k-1}=k
\big\}.
\end{eqnarray*}
On the other hand, consider the representation $S(1)=I_1({\bf 1}_{[0,1]})$.
As above, iterative applications of the product formula (\ref{prodfree}) leads to
\[
S(1)^k=I_1({\bf 1}_{[0,1]})^k=\sum_{(r_1,\ldots,r_{k-1})\in \widetilde{A}_k}
I_{k-2r_1-\ldots-2r_{k-1}}\big(
{\bf 1}_{[0,1]}
\overset{r_1}{\frown}
\ldots \overset{r_{k-1}}{\frown}{\bf 1}_{[0,1]}
\big),
\]
where
\begin{eqnarray*}
\widetilde{A}_k&=&\big\{
(r_1,\ldots,r_{k-1})\in \{0,1\}^{k-1}:\,r_2\leq 2-2r_1,\,\,r_3\leq 3-2r_1-2r_2,\ldots,\\
&&\hskip7cm r_{k-1}\leq k-1-2r_1-\ldots-2r_{k-2}
\big\}.
\end{eqnarray*}
By taking the expectation, we deduce that
\[
E[S(1)^k]=\sum_{(r_1,\ldots,r_{k-1})\in \widetilde{C}_k}
{\bf 1}_{[0,1]}\overset{r_1}{\frown}\ldots \overset{r_{k-1}}{\frown}{\bf 1}_{[0,1]}
=\sum_{(r_1,\ldots,r_{k-1})\in \widetilde{C}_k} 1
= \# \widetilde{C}_k = \# C_k.
\]
\fin
\begin{rem}
{\rm
When $k$ is even, it is well-known that $E[S(1)^k]$ is given
by ${\rm Cat}_{k/2}$, the Catalan number of order $k/2$.
There is many combinatorial ways to define this number. One of them is to see
it at the number of paths
in the lattice $\mathbb{Z}^2$ which start at $(0,0)$, end at $(k,0)$, make
steps of the form $(1,1)$ or $(1,-1)$, and never lies below the $x$-axis, i.e., all
their points are of the form $(i,j)$ with $j\geq 0$.

Let the notation of the proof of Lemma \ref{lm2} prevail.
Set $s_i=1-2r_i$. Then
\begin{eqnarray*}
\widetilde{C}_k&\overset{\rm bij.}{\equiv}&\left\{
(s_1,\ldots,s_{k-1})\in \{-1,1\}^{k-1}:\,1+s_1\geq\frac12(1-s_2),\,1+s_1+s_2\geq \frac12(1-s_3),\right.\\
&&\hskip3.5cm \ldots,1+s_1+\ldots+s_{k-2}\geq \frac12(1-s_{k-1}),\,1+s_1+\ldots+s_{k-1}=0
\bigg\}.
\end{eqnarray*}
It turns out that the set of conditions
\begin{equation}\label{cond1}
\left\{
\begin{array}{lll}
s_j\in\{-1,1\},\quad j=1,\ldots,k-1\\
1+s_1+\ldots+s_j\geq \frac12(1-s_{j+1}),\quad j=1,\ldots,k-2\\
1+s_1+\ldots+s_{k-1}=0,
\end{array}
\right.
\end{equation}
is equivalent to
\begin{equation}\label{cond2}
\left\{
\begin{array}{lll}
s_j\in\{-1,1\},\quad j=1,\ldots,k-1\\
1+s_1+\ldots+s_j\geq 0,\quad j=1,\ldots,k-2\\
1+s_1+\ldots+s_{k-1}=0.
\end{array}
\right.
\end{equation}
Indeed, it is clear that (\ref{cond1}) implies (\ref{cond2}).
Conversely, suppose that (\ref{cond2}) is in order, and let $j\in\{1,\ldots,k-2\}$.
Because  $\frac12(1-s_{j+1})\leq 1$, one has that $1+s_1+\ldots+s_j\geq \frac12(1-s_{j+1})$
when $1+s_1+\ldots+s_j\geq 1$.
If $1+s_1+\ldots+s_j=0$ then, because $1+s_1+\ldots+s_{j+1}\geq 0$ (even if $j=k-2$),
one has $s_{j+1}=1$, implying in turn $1+s_1+\ldots+s_j\geq \frac12(1-s_{j+1})=0$.
Thus
\begin{eqnarray*}
\widetilde{C}_k&\overset{\rm bij.}{\equiv}&\left\{
(s_1,\ldots,s_{k-1})\in \{-1,1\}^{k-1}:\,1+s_1 \geq 0,\,1+s_1+s_2\geq 0,\right.\\
&&\hskip3.5cm \ldots,1+s_1+\ldots+s_{k-2}\geq 0,\,1+s_1+\ldots+s_{k-1}=0
\bigg\},
\end{eqnarray*}
and we recover the result of Lemma \ref{lm2} when $k$ is even. (The case where $k$ is odd is
trivial.)
}
\end{rem}

\begin{lemme}\label{lm3}
We have $f_n\overset{r_1}{\frown}\ldots\overset{r_{k-1}}{\frown}f_n=1$
for all $k\geq 3$ and all $(r_1,\ldots,r_{k-1})\in C_k$.
\end{lemme}
{\it Proof}.
It is evident, using the identities
$f_n\overset{0}{\frown} f_n = f_n\otimes f_n$ and \[
f_n\overset{p}{\frown} f_n =
\int_{[0,1]^p}f_n(t_1,\ldots,t_p)f_n(t_p,\ldots,t_1)dt_1\ldots dt_p=\|f_n\|^2_{L^2([0,1]^p)}=1.\]
\fin
\begin{lemme}\label{lm4}
As $n\to\infty$, assume that $\|f_n\overset{r}{\frown} f_n\|_{L^2([0,1]^{2p-2r})}\to 0$ for all $r=1,\ldots,p-1$.
Then, as $n\to\infty$ we have $
f_n\overset{r_1}{\frown}\ldots\overset{r_{k-1}}{\frown}f_n
\to 0$
for all $k\geq 3$ and all $(r_1,\ldots,r_{k-1})\in E_k$.
\end{lemme}
{\it Proof}. Fix $(r_1,\ldots,r_{k-1})\in E_k$, and
let $j\in\{1,\ldots,k-1\}$ be the smallest integer such that $r_j\in\{1,\ldots,p-1\}$.
Recall that $f_n\overset{0}{\frown} f_n = f_n\otimes f_n$.
Then
\begin{eqnarray*}
&&\big| f_n\overset{r_1}{\frown}\ldots\overset{r_{k-1}}{\frown}f_n \big| \\
&=&\big|
f_n\overset{r_1}{\frown}\ldots   \overset{r_{j-1}}{\frown}  f_n \overset{r_{j}}{\frown} f_n
\overset{r_{j+1}}{\frown} \ldots \overset{r_{k-1}}{\frown} f_n\big|\\
&=&\big|
(f_n\otimes\ldots \otimes f_n) \overset{r_{j}}{\frown} f_n
\overset{r_{j+1}}{\frown} \ldots \overset{r_{k-1}}{\frown} f_n\big|
\quad
\mbox{(using $f_n\overset{p}{\frown} f_n=1$)}\\
&\leq&
\|(f_n\otimes\ldots \otimes f_n)\otimes (f_n  \overset{r_{j}}{\frown} f_n)\|_{L^2([0,1]^q)}
\|f_n\|^{k-j-1}_{L^2([0,1]^p)}\quad
\mbox{(by Cauchy-Schwarz, for a certain $q$)}\\
&=&\|f_n  \overset{r_{j}}{\frown} f_n\|\quad\mbox{(because $\|f_n\|_{L^2([0,1]^p)}^2=1$)}\\
&\longrightarrow& 0\quad\mbox{as $n\to\infty$}.
\end{eqnarray*}
\fin

\section{Proof of Theorem  \ref{nunugio}}
We follow the same route as in the proof of Theorem \ref{knps}, that is, we utilize the method of moments.
(It is well-known that the $N(0,1)$ law is uniquely determined by its moments.)
Let the notation and assumptions of Theorem \ref{nunugio} prevail.
Without loss of generality, we may assume that $E[F_n^2]=1$ for all $n$
(instead of $E[F_n^2]\to 1$ as $n\to\infty$).
Moreover, observe that
$p!\|f_n\|^2_{L^2([0,1]^p)}=E[F_n^2]=1$.

Fix an integer $k\geq 3$. Iterative applications of the product formula (\ref{product}) leads to
\begin{eqnarray}\label{kthmomentwithoutexpectation-bis}
F_n^k=I_p(f_n)^k&=&\sum_{(r_1,\ldots,r_{k-1})\in A_k}
I_{kp-2r_1-\ldots-2r_{k-1}}\big(
f_n\widetilde{\otimes}_{r_1}\ldots\widetilde{\otimes}_{r_{k-1}}f_n
\big)\\
&&\hskip5cm\times\prod_{j=1}^{k-1}r_j!\binom{p}{r_j}\binom{jp-2r_1-\ldots-2r_{j-1}}{r_j},\notag
\end{eqnarray}
where
\begin{eqnarray*}
A_k&=&\big\{
(r_1,\ldots,r_{k-1})\in \{0,1,\ldots,p\}^{k-1}:\,r_2\leq 2p-2r_1,\,\,r_3\leq 3p-2r_1-2r_2,\ldots,\\
&&\hskip7cm r_{k-1}\leq(k-1)p-2r_1-\ldots-2r_{k-2}
\big\}.
\end{eqnarray*}
In order to simplify the exposition, note that we have removed all the brackets in the writing of
$f_n
\widetilde{\otimes}_{r_1}
\ldots
\widetilde{\otimes}_{r_{k-1}}
f_n$. We use the implicit convention
that
these quantities are always  defined iteratively from the left to the right. For instance,
$f_n\widetilde{\otimes}_{r_1}f_n \widetilde{\otimes}_{r_2}f_n \widetilde{\otimes}_{r_3}f_n$
stands for $((f_n\widetilde{\otimes}_{r_1}f_n) \widetilde{\otimes}_{r_2}f_n) \widetilde{\otimes}_{r_3}f_n$.

By taking the expectation in (\ref{kthmomentwithoutexpectation-bis}), we deduce that
\begin{equation}\label{kthmoment-bis}
E[F_n^k]=\sum_{(r_1,\ldots,r_{k-1})\in B_k}
f_n
\widetilde{\otimes}_{r_1}
\ldots
\widetilde{\otimes}_{r_{k-1}}
f_n\times\prod_{j=1}^{k-1}r_j!\binom{p}{r_j}\binom{jp-2r_1- \ldots-2r_{j-1}}{r_j},
\end{equation}
with $B_k=\big\{(r_1,\ldots,r_{k-1})\in A_k:\,2r_1+\ldots+2r_{k-1}=kp\big\}$.
Combining (\ref{kthmoment-bis}) with the crude bound (consequence of Cauchy-Schwarz)
\[
\|f_n \widetilde{\otimes}_r f_n\|_{L^2([0,1]^{2p-2r})}
\leq \|f_n\|^2_{L^2([0,1]^p)}=1/p!\leq 1,\]
we have that $E[F_n^k]\leq \# B_k$, that is, for every $k$ the $k$th
moment of $F_n$ is uniformly bounded.

Assume that $(i)$ is in order.
Because of the uniform boundedness of the moments, standard arguments implies that
$E[F_n^4]\to E[B(1)^4]$.
Conversely, assume that $(ii)$ is in order and let us prove that, for all $k\geq 1$,
\begin{equation}\label{tobeshown-bis}
E[F_n^k]\to E[B(1)^k]\quad\mbox{as $n\to\infty$.}
\end{equation}
The cases $k=1$ and $k=2$ being immediate, assume that $k\geq 3$ is given.
We decompose $B_k$ as $C_k\cup E_k$, with $C_k=B_k\cap\{0,p\}^{k-1}$ and $E_k=B_k\setminus C_k$.
We have
\begin{eqnarray}
E[F_n^k]&=& \sum_{(r_1,\ldots,r_{k-1})\in C_k}
f_n\widetilde{\otimes}_{r_1} \ldots\widetilde{\otimes}_{r_{k-1}}  f_n\times
\prod_{j=1}^{k-1}r_j!\binom{jp-2r_1-\ldots-2r_{j-1}}{r_j}\label{doublesum-bis}
\\
&&+\sum_{(r_1,\ldots,r_{k-1})\in E_k}
f_n\widetilde{\otimes}_{r_1}  \ldots\widetilde{\otimes}_{r_{k-1}}  f_n
\times \prod_{j=1}^{k-1}r_j!\binom{p}{r_j}\binom{jp-2r_1-\ldots-2r_{j-1}}{r_j}.\notag
\end{eqnarray}
By Lemma \ref{lm1-bis} together with assumption $(ii)$, we have
that $\|f_n\otimes_r f_n\|_{L^2([0,1]^{2p-2r})}$
(as well as $\|f_n\widetilde{\otimes}_r f_n\|_{L^2([0,1]^{2p-2r})}$) tends to zero for any $r=1,\ldots,p-1$.
Lemmas \ref{lm0} and \ref{lm2-bis} imply together that
the first sum in (\ref{doublesum-bis}) converges to $E[B(1)^k]$, whereas
the second sum converges to zero by Lemma \ref{lm4-bis}.
Thus, (\ref{tobeshown-bis}) is in order, and the proof of the theorem is concluded.\fin

\begin{lemme}\label{lm1-bis}
We have \[
E[F_n^4]=3+\sum_{r=1}^{p-1} \binom{p}{r}^2
\left[(p!)^2
\|f_n\otimes_r f_n\|^2_{L^2([0,1]^{2p-2r})}
+(r!)^2\binom{p}{r}^2(2p-2r)!
\|f_n \widetilde{\otimes}_r\, f_n\|^2_{L^2([0,1]^{2p-2r})}
\right].
\]
\end{lemme}
{\it Proof} (following \cite{nualart-peccati}).
Let $\pi\in\mathfrak{S}_{2p}$.
If $r\in\{0,\ldots,p\}$ denotes the cardinality of $\{\pi(1),\ldots,\pi(p)\}\cap\{1,\ldots,p\}$
then it is readily checked that $r$ is also the cardinality of
$\{\pi(p+1),\ldots,\pi(2p)\}\cap\{p+1,\ldots,2p\}$ and that
\begin{eqnarray}
&&\int_{[0,1]^{2p}}f_n(t_1,\ldots,t_p)f_n(t_{\pi(1)},\ldots,t_{\pi(p)})f_n(t_{p+1},\ldots,t_{2p})
f_n(t_{\pi(p+1)},\ldots,t_{\pi(2p)})dt_1\ldots dt_{2p}\notag\\
&=&\int_{[0,1]^{2p-2r}}f_n\otimes_r f_n(x_1,\ldots,x_{2p-2r})^2dx_1\ldots dx_{2p-2r} =
\|f_n\otimes_r f_n\|^2_{L^2([0,1]^{2p-2r})}.\label{ctr}
\end{eqnarray}
Moreover, for any fixed $r\in\{0,\ldots,p\}$, there are $\binom{p}{r}^2(p!)^2$
permutations $\pi\in\mathfrak{S}_{2p}$ such that
$\#\{\pi(1),\ldots,\pi(p)\}\cap\{1,\ldots,p\}=r$.
(Indeed, such a permutation is completely determined by the choice of: $(a)$ $r$ distinct
elements $x_1,\ldots,x_r$ of $\{1,\ldots,p\}$; $(b)$ $p-r$ distinct elements $x_{r+1},\ldots,x_p$
of $\{p+1,\ldots,2p\}$; $(c)$ a bijection between $\{1,\ldots,p\}$ and $\{x_1,\ldots,x_p\}$;
$(d)$ a bijection between $\{p+1,\ldots,2p\}$ and $\{1,\ldots,2p\}\setminus \{x_1,\ldots,x_p\}$.)
Now, recall from (\ref{norm}) that the symmetrization of $f_n\otimes f_n$ is given by
\[
f_n\widetilde{\otimes} f_n(t_1,\ldots,t_{2p}) = \frac{1}{(2p)!}
\sum_{\pi\in\mathfrak{S}_{2p}} f_n(t_{\pi(1)},\ldots,t_{\pi(p)})
f_n(t_{\pi(p+1)},\ldots,t_{\pi(2p)}).
\]
Therefore,
\begin{eqnarray*}
\|f_n\widetilde{\otimes} f_n\|^2_{L^2([0,1]^{2p})}&=&
\frac{1}{(2p)!^2}
\sum_{\pi,\pi'\in\mathfrak{S}_{2p}}
\int_{[0,1]^{2p}}
f_n(t_{\pi(1)},\ldots,t_{\pi(p)})f_n(t_{\pi(p+1)},\ldots,t_{\pi(2p)})\\
&&\hskip3cm\times
f_n(t_{\pi'(1)},\ldots,t_{\pi'(p)})f_n(t_{\pi'(p+1)},\ldots,t_{\pi'(2p)})
dt_1\ldots dt_{2p}\\
&=&
\frac{1}{(2p)!}
\sum_{\pi\in\mathfrak{S}_{2p}}
\int_{[0,1]^{2p}}
f_n(t_{1},\ldots,t_{p})f_n(t_{p+1},\ldots,t_{2p})\\
&&\hskip3cm\times
f_n(t_{\pi(1)},\ldots,t_{\pi(p)})f_n(t_{\pi(p+1)},\ldots,t_{\pi(2p)})
dt_1\ldots dt_{2p}\\
&=&\frac{1}{(2p)!}\sum_{r=0}^p
\sum_{\substack{\pi\in\mathfrak{S}_{2p}\\
\{\pi(1),\ldots,\pi(p)\}\cap\{1,\ldots,p\}=r
}}
\int_{[0,1]^{2p}}
f_n(t_{1},\ldots,t_{p})f_n(t_{p+1},\ldots,t_{2p})\\
&&\hskip3cm\times
f_n(t_{\pi(1)},\ldots,t_{\pi(p)})f_n(t_{\pi(p+1)},\ldots,t_{\pi(2p)})
dt_1\ldots dt_{2p}.
\end{eqnarray*}
Hence, using (\ref{ctr}), we deduce that
\begin{eqnarray}
(2p)!\|f_n\widetilde{\otimes} f_n\|^2_{L^2([0,1]^{2p})}
&=&2(p!)^2\|f_n\|_{L^2([0,1]^p)}^4+(p!)^2\sum_{r=1}^{p-1}
\binom{p}{r}^2
\|f_n\otimes_r f_n\|^2_{L^2([0,1]^{2p-2r})}
\notag\\
&=&2+(p!)^2\sum_{r=1}^{p-1}
\binom{p}{r}^2\|f_n\otimes_r f_n\|^2_{L^2([0,1]^{2p-2r})}.\label{beautyformula}
\end{eqnarray}
The product formula (\ref{product}) leads to
$
F_n^2= \sum_{r=0}^p r!\binom{p}{r}^2 I_{2p-2r}(f_n\widetilde{\otimes}_r\, f_n).
$
Using (\ref{isom})-(\ref{isom-dif}), we infer
\begin{eqnarray*}
E[F_n^4] &=& \sum_{r=0}^{p} (r!)^2\binom{p}{r}^4 (2p-2r)!
\|f_n\widetilde{\otimes}_r f_n\|^2_{L^2([0,1]^{2p-2r})}\\
&=&(2p)! \|f_n\widetilde{\otimes} f_n\|^2_{L^2([0,1]^{2p})}
+1
+\sum_{r=1}^{p-1} (r!)^2\binom{p}{r}^4 (2p-2r)!
\|f_n\widetilde{\otimes}_r f_n\|^2_{L^2([0,1]^{2p-2r})}.
\end{eqnarray*}
By inserting (\ref{beautyformula}) in the previous identity, we get the desired result.
\fin
\begin{lemme}\label{lm0}
As $n\to\infty$, assume that
\begin{equation}\label{labellecontr}
\|f_n\otimes_r f_n\|_{L^2([0,1]^{2p-2r})}\to 0,\quad r=1,\ldots,p-1.
\end{equation}
Then, for all $k\geq 3$ and all $(r_1,\ldots,r_{k-1})\in C_k$, we have
\[
f_n\widetilde{\otimes}_{r_1} \ldots\widetilde{\otimes}_{r_{k-1}}  f_n\to
\prod_{j=1}^{k-1} \frac{\binom{j-2r_1/p-\ldots-2r_{j-1}/p}{r_j/p}}
{(r_j)!\binom{jp-2r_1-\ldots-2r_{j-1}}{r_j}}
\quad\mbox{as $n\to\infty$}.
\]
\end{lemme}
{\it Proof}. In all the proof, for sake of conciseness we write
$f_n^{\widetilde{\otimes}d}$ instead of $\overbrace{f_n\widetilde{\otimes}\ldots\widetilde{\otimes}f_n}^{\mbox{$d$ times}}$.
(Here, ``$d$ times''
just means that $f_n$ appears $d$ times in the expression.)
It is readily checked that $f_n^{\widetilde{\otimes}d}=\widetilde{f_n^{\otimes d}}$
so that, according to (\ref{norm}),
\begin{eqnarray*}
f_n^{\widetilde{\otimes}d}\otimes_p f_n(t_1,\ldots,t_{dp-p})
&=&\frac{1}{(dp)!}\sum_{\pi\in\mathfrak{S}_{dp}}\int_{[0,1]^{p}}
f_n(t_{\pi(1)},\ldots,t_{\pi(d)})\ldots f_n(t_{\pi(dp-p+1)},\ldots,t_{\pi(dp)})\\
&&\hskip3cm\times f_n(t_{dp-p+1},\ldots,t_{dp})dt_{dp-d+1}\ldots dt_{dp}.
\end{eqnarray*}
Let $\pi\in\mathfrak{S}_{dp}$.
When
$\{\pi(jp-p+1),\ldots,\pi(jp)\}\neq \{ dp-p+1,\ldots,dp\}$
for all $j=1,\ldots,d$,
it is readily checked, using
(\ref{labellecontr}) as well as Cauchy-Schwarz, that
 the function
\begin{eqnarray*}
(t_1,\ldots,t_{dp-p})&\mapsto& \int_{[0,1]^{p}}
f_n(t_{\pi(1)},\ldots,t_{\pi(d)})\ldots f_n(t_{\pi(dp-p+1)},\ldots,t_{\pi(dp)}) \\
&&\hskip5cm\times
f_n(t_{dp-p+1},\ldots,t_{dp})dt_{dp-d+1}\ldots dt_{dp}
\end{eqnarray*}
tends to zero in $L^2([0,1]^{dp-p})$.
Let $\mathfrak{A}_{dp}$ be the set of permutations $\pi\in\mathfrak{S}_{dp}$
for which there exists (at least one) $j\in \{1,\ldots,d\}$ such that
$\{\pi(jp-p+1),\ldots,\pi(jp)\}=\{ dp-p+1,\ldots,dp\}$.
We then have
\begin{eqnarray*}
f_n^{\widetilde{\otimes}d}\otimes_p f_n(t_1,\ldots,t_{dp-p})
&\approx&\frac{1}{(dp)!}\sum_{\pi\in\mathfrak{A}_{dp}}\int_{[0,1]^{p}}
f_n(t_{\pi(1)},\ldots,t_{\pi(d)})\ldots f_n(t_{\pi(dp-p+1)},\ldots,t_{\pi(dp)})\notag\\
&&\hskip3cm\times f_n(t_{dp-p+1},\ldots,t_{dp})dt_{dp-d+1}\ldots dt_{dp},
\end{eqnarray*}
where, here and in the rest of the proof, we use the notation $h_n\approx g_n$ (for $h_n$ and $g_n$ two functions
of, say, $q$ arguments) to mean that $h_n-g_n$
tends to zero in $L^2([0,1]^{q})$.
Because a permutation $\pi$ of $\mathfrak{A}_{dp}$
is completely characterized by the choice of the smallest index
$j$ for which $\{\pi(jp-p+1),\ldots,\pi(jp)\}=\{ dp-p+1,\ldots,dp\}$ as well as two permutations $\tau\in\mathfrak{S}_{p}$
and $\sigma\in\mathfrak{S}_{pd-p}$, and using moreover that $f_n\otimes_p f_n = \|f_n\|^2_{L^2([0,1]^p)}=\frac{1}{p!}$
and that $f_n$ is symmetric, we deduce that
\begin{eqnarray}
f_n^{\widetilde{\otimes}d}\otimes_p f_n(t_1,\ldots,t_{dp-p})
&\approx&
\frac{d}{(dp)!} \sum_{\sigma\in\mathfrak{S}_{dp-p}}
f_n(t_{\sigma(1)},\ldots,t_{\sigma(d)})\ldots f_n(t_{\sigma(dp-2p+1)},\ldots,t_{\sigma(dp-p)}) \notag\\
&\approx& \frac{d}{p!\binom{dp}{p}}\widetilde{f_n^{\otimes (d-1)}}(t_1,\ldots,t_{dp-p})
=\frac{d}{p!\binom{dp}{p}}f_n^{\widetilde{\otimes} (d-1)}(t_1,\ldots,t_{dp-p}).\notag\\
\label{crucru}
\end{eqnarray}
Because the right-hand side of (\ref{crucru}) is a symmetric function, we eventually
get that
\[
f_n^{\widetilde{\otimes}d}\widetilde{\otimes}_p f_n \approx
\frac{d}{p!\binom{dp}{p}}f_n^{\widetilde{\otimes} (d-1)},
\]
with the convention that $f_n^{\widetilde{\otimes}0}=1$.
On the other hand, we have
$f_n^{\widetilde{\otimes}d}\widetilde{\otimes}_0 f_n=f_n^{\widetilde{\otimes}d}\widetilde{\otimes} f_n = f_n^{\widetilde{\otimes}(d+1)}$ by the very
definition of $f_n^{\widetilde{\otimes}d}$.
We can summarize these two last identities by writing that, for any $r\in\{0,p\}$,
\begin{equation}\label{cru5}
f_n^{\widetilde{\otimes}d}\widetilde{\otimes}_r f_n \approx
\frac{\binom{d}{r/p}}{r!\binom{dp}{r}}\,
f_n^{\widetilde{\otimes}(d+1-2r/p)}.
\end{equation}
Now, let $k\geq 3$ and $(r_1,\ldots,r_{k-1})\in C_k$.
Thanks to (\ref{cru5}), we have
$
f_n\widetilde{\otimes}_{r_1} f_n =
\frac{\binom{1}{r_1/p}}{(r_1)!\binom{p}{r_1}}\,
f_n^{\widetilde{\otimes}(2-2r_1/p)},
$
\begin{eqnarray*}
f_n\widetilde{\otimes}_{r_1} f_n \widetilde{\otimes}_{r_2} f_n \approx
\frac{\binom{1}{r_1/p}\binom{2-2r_1/p}{r_2/p}}{(r_1)!\binom{p}{r_1}(r_2)!\binom{2p-2r_1}{r_2}}
f_n^{\widetilde{\otimes}(3-2r_1/p-2r_2/p)},
\end{eqnarray*}
and so on. Iterating this procedure leads eventually to
\begin{equation}\label{crul}
f_n\widetilde{\otimes}_{r_1}\ldots\widetilde{\otimes}_{r_{k-1}}f_n
\approx
\prod_{j=1}^{k-1} \frac{\binom{j-2r_1/p-\ldots-2r_{j-1}/p}{r_j/p}}
{(r_j)!\binom{jp-2r_1-\ldots-2r_{j-1}}{r_j}},
\end{equation}
which is exactly the desired formula. The proof of the lemma is done.\fin
\begin{lemme}\label{lm2-bis}
For all $k\geq 3$, we have
\[
 E[B(1)^k]=\sum_{(r_1,\ldots,r_{k-1})\in C_k}
\prod_{j=1}^{k-1}\binom{j-2r_1/p-\ldots-2r_{j-1}/p}{r_j/p}.
\]
\end{lemme}
{\it Proof}.
The identity is clear when $k$ is an odd integer, because $C_k=\emptyset$ in this case.
Assume now that $k$ is even.
Consider the representation $B(1)=I_1({\bf 1}_{[0,1]})$.
Iterative applications of the product formula (\ref{product}) leads to
\begin{eqnarray*}
B(1)^k&=&I_1({\bf 1}_{[0,1]})^k=\sum_{(r_1,\ldots,r_{k-1})\in \widetilde{A}_k}
I_{k-2r_1-\ldots-2r_{k-1}}\big(
{\bf 1}_{[0,1]}
\widetilde{\otimes}_{r_1}
\ldots
\widetilde{\otimes}_{r_{k-1}}{\bf 1}_{[0,1]}
\big)\\
&&\hskip8cm \times \prod_{j=1}^{k-1}\binom{j-2r_1-\ldots-2r_{j-1}}{r_j},
\end{eqnarray*}
where
\begin{eqnarray*}
\widetilde{A}_k&=&\big\{
(r_1,\ldots,r_{k-1})\in \{0,1\}^{k-1}:\,r_2\leq 2-2r_1,\,\,r_3\leq 3-2r_1-2r_2,\ldots,\\
&&\hskip7cm r_{k-1}\leq k-1-2r_1-\ldots-2r_{k-2}
\big\}.
\end{eqnarray*}
By taking the expectation, we deduce that
\[
E[B(1)^k]=\sum_{(r_1,\ldots,r_{k-1})\in \widetilde{C}_k}
{\bf 1}_{[0,1]}
\widetilde{\otimes}_{r_1}
\ldots
\widetilde{\otimes}_{r_{k-1}}
{\bf 1}_{[0,1]}\times\prod_{j=1}^{k-1}\binom{j-2r_1-\ldots-2r_{j-1}}{r_j},
\]
with
\begin{eqnarray*}
\widetilde{C}_k&=&\big\{
(r_1,\ldots,r_{k-1})\in \{0,1\}^{k-1}:\,r_2\leq 2-2r_1,\,\,r_3\leq 3-2r_1-2r_2,\ldots,\\
&&\hskip3.5cm r_{k-1}\leq k-1-2r_1-\ldots-2r_{k-2},\,2r_1+\ldots+2r_{k-1}=k
\big\}.
\end{eqnarray*}
It is readily checked that
$
{\bf 1}_{[0,1]}
\widetilde{\otimes}_{r_1}
\ldots
\widetilde{\otimes}_{r_{k-1}}
{\bf 1}_{[0,1]} = {\bf 1}_{[0,1]}
\otimes_{r_1}
\ldots
\otimes_{r_{k-1}}
{\bf 1}_{[0,1]} =1$
for all $(r_1,\ldots,r_{k-1})\in \widetilde{C}_k$.
Hence
\begin{eqnarray*}
E[B(1)^k]&=&\sum_{(r_1,\ldots,r_{k-1})\in \widetilde{C}_k}
\prod_{j=1}^{k-1}\binom{j-2r_1-\ldots-2r_{j-1}}{r_j}\\
&=&\sum_{(r_1,\ldots,r_{k-1})\in C_k}
\prod_{j=1}^{k-1}\binom{j-2r_1/p-\ldots-2r_{j-1}/p}{r_j/p},
\end{eqnarray*}
which is the desired conclusion.
\fin
\begin{lemme}\label{lm4-bis}
As $n\to\infty$, assume that $\|f_n\widetilde{\otimes}_r f_n\|_{L^2([0,1]^{2p-2r})}\to 0$ for all $r=1,\ldots,p-1$.
Then, as $n\to\infty$ we have $
f_n
\widetilde{\otimes}_{r_1}
\ldots
\widetilde{\otimes}_{r_{k-1}}
f_n
\to 0$
for all $k\geq 3$ and all $(r_1,\ldots,r_{k-1})\in E_k$.
\end{lemme}
{\it Proof}. Fix $k\geq 3$ and $(r_1,\ldots,r_{k-1})\in E_k$, and
let $j\in\{1,\ldots,k-1\}$ be the smallest integer such that $r_j\in\{1,\ldots,p-1\}$.
As in the proof of Lemma \ref{lm0}, when $h_n$ and $g_n$ are functions of $q$ arguments let us write
$h_n\approx g_n$ to indicate that $h_n-g_n$ tends to zero in $L^2([0,1]^q)$.
Recall from (\ref{cru5}) that $f_n^{\widetilde{\otimes} d}\widetilde{\otimes}_p f_n \approx \frac{d}{p!\binom{dp}{p}}f_n^{\widetilde{\otimes}(d-1)}$.
Then
\begin{eqnarray*}
&&\big| f_n\widetilde{\otimes}_{r_{1}}\ldots\widetilde{\otimes}_{r_{k-1}}f_n \big| \\
&=&\big|
f_n\widetilde{\otimes}_{r_{1}}\ldots   \widetilde{\otimes}_{r_{j-1}}  f_n \widetilde{\otimes}_{r_{j}} f_n
\widetilde{\otimes}_{r_{j+1}} \ldots \widetilde{\otimes}_{r_{k-1}}f_n\big|\\
&\approx &c\big|
(f_n\widetilde{\otimes}\ldots \widetilde{\otimes}f_n) \widetilde{\otimes}_{r_{j}} f_n
\widetilde{\otimes}_{r_{j+1}} \ldots \widetilde{\otimes}_{r_{k-1}} f_n\big|
\quad
\mbox{(for some constant $c>0$ independent of $n$)}\\
&\leq&c
\|(f_n\widetilde{\otimes}\ldots \widetilde{\otimes} f_n)\widetilde{\otimes} (f_n  \widetilde{\otimes}_{r_{j}} f_n)\|_{L^2([0,1]^q)}
\|f_n\|^{k-j-1}_{L^2([0,1]^p)}\quad
\mbox{(by Cauchy-Schwarz, for a certain $q$)}\\
&\leq &c\|f_n  \widetilde{\otimes}_{r_{j}} f_n\|
\quad\mbox{(because $\|f_n\|_{L^2([0,1]^p)}^2=\frac{1}{p!}\leq 1$)}\\
&\longrightarrow& 0\quad\mbox{as $n\to\infty$}.
\end{eqnarray*}
\fin

\bigskip
\noindent
{\bf Acknowledgement}. I thank one anonymous referee for his/her thorough reading and insightful comments.

\end{document}